\documentclass[11pt]{amsart}

\usepackage{latexsym}
\usepackage{amsfonts}
\usepackage{amsthm}
\usepackage{graphicx}
\usepackage{amsmath}
\usepackage{amssymb}
\usepackage{color}
\newtheorem{problem}{Problem}

\newcommand{\gp}[1]{{\left\langle #1 \right\rangle}}
\begin{document}

\title{Search and witness problems in group theory}

\author{Vladimir Shpilrain}
\address{Department of Mathematics, The City  College  of New York, New York,
NY 10031} \email{shpil@groups.sci.ccny.cuny.edu}
\thanks{Partially supported by the NSF grant DMS-0914778. }

\begin{abstract}
Decision problems are problems of the following nature: given a
property $\mathcal{P}$ and an object $\mathcal{O}$, find out whether
or not the object $\mathcal{O}$ has the property $\mathcal{P}$. On
the other hand, witness problems are: given a property $\mathcal{P}$
and an  object $\mathcal{O}$ with the property $\mathcal{P}$, find a
proof of the fact that $\mathcal{O}$ indeed has the property
$\mathcal{P}$.

On the third hand(?!), search problems are  of the following nature:
given a property $\mathcal{P}$  and an  object $\mathcal{O}$ with
the property $\mathcal{P}$, find something ``material" establishing
the property $\mathcal{P}$; for example, given two conjugate
elements of a group, find a  conjugator. In this survey our focus is
on various search problems in group theory, including the word
search problem, the subgroup membership search problem, the
conjugacy search problem, and others.

\end{abstract}

\maketitle

\begin{center}

{\sc \small \it  To Alfred Lvovich Shmelkin with deepest
appreciation}

\end{center}

\bigskip

\section{Introduction}

Decision problems are problems of the following nature: given a
property $\mathcal{P}$ and an object $\mathcal{O}$, find out whether
or not the object $\mathcal{O}$ has the property $\mathcal{P}$. On
the other hand, {\it  search problems} are  of the following nature:
given a property $\mathcal{P}$  and an  object $\mathcal{O}$ with
the property $\mathcal{P}$, find something ``material"  establishing
the property $\mathcal{P}$; for example, given two conjugate
elements of a group, find a  conjugator. A weaker version of a
search problem is sometimes called a {\it witness problem}: given a
property $\mathcal{P}$ and an  object $\mathcal{O}$ with the
property $\mathcal{P}$, find a proof of the fact that $\mathcal{O}$
indeed has the property $\mathcal{P}$.

Search  and   witness   problems represent a substantial shift of
paradigm from decision problems, and in fact, studying witness and
search problems often gives rise to new research avenues in
mathematics, very different from those prompted by addressing  the
corresponding decision problems. To give just a couple of examples
from different areas of mathematics, we can mention (1) the
isoperimetric function that can be used to measure the complexity of
a proof that a given word is trivial in a given group; (2)
Reidemeister moves that can be used to measure the complexity of a
proof that two given knot diagrams are those of two isotopic knots;
(3) elementary row (or column) operations on a matrix over a field
that can be used to measure the complexity of a proof that a given
square matrix  is invertible. With respect to the last example we
note that, although a more straightforward proof would be producing
the inverse matrix (this would solve the relevant search problem),
the proof by elementary row (or column) operations provides a useful
{\it stratification} of the relevant witness problem, which allows
one to allocate a witness problem to one of the established
complexity classes (e.g. ${\mathbf P}$ or ${\mathbf {NP}}$) by
converting it to a decision problem; in this particular example the
latter would be asking whether or not a given matrix is a product of
at most $k$ elementary matrices.

The main objective of this survey is to discuss various search
problems in group theory. We note that decision problems  in group
theory have been studied for over 100 years now, since Dehn put
forward, in the beginning of the 20th century, the three famous
decision problems now often referred to as {\it Dehn's problems}:
the word problem, the conjugacy problem, and the isomorphism
problem. Later, some of these problems were generalized, and many
other decision problems were raised; we refer to   \cite{KS-survey}
or \cite{cfm-survey} for a survey.

On the other hand, search problems in group theory and their
complexity started to attract attention relatively recently.
Complexity of the word search problem in a finitely presented group
is reflected by isoperimetric and isodiametric functions of a finite
presentation of this group, as introduced in \cite{Gromov} and
\cite{Gersten} in 1985--1991. More recently, complexity of the
conjugacy search problem has got a lot of attention, after a seminal
paper \cite{AAG} offered a cryptographic key exchange protocol that
relied in its security on the complexity of the conjugacy search
problem in braid groups.

Later on, there were other proposals of cryptographic primitives
that relied in their security on the complexity of other search
problems (see \cite{MSU_book} for a comprehensive survey), including
the word search problem, the subgroup membership search problem
\cite{SZ}, the decomposition search problem \cite{SU3}, etc. This
has boosted interest in studying  various search problems in groups,
and it is the purpose of the present survey to expose at least some
of the directions of this research.

\section{Decision and  search problems in group theory}
\label{alg_problems}

As we have pointed out in the Introduction, algorithmic problems
considered in group theory  are of three different kinds:

\begin{itemize}

\item {\it Decision problems}  are
problems of the following nature: given a property $\mathcal{P}$ and
an object $\mathcal{O}$, find out whether or not the  object
$\mathcal{O}$ has the property $\mathcal{P}$.

\item {\it Witness problems} are  of the
following nature: given a property $\mathcal{P}$  and an  object
$\mathcal{O}$ with the property $\mathcal{P}$, find a proof (a
``witness") of the fact that $\mathcal{O}$  has the property
$\mathcal{P}$. Such a proof does not  necessarily have to produce
anything ``material"; for example, we mentioned in the Introduction
that one of the ways to prove invertibility of a matrix over a field
is reducing it by elementary row or column operations to the
identity matrix. This way does not by itself produce the inverse of
a given matrix, although, of course, upon some little extra effort
it will.

\item {\it Search problems} are, typically, a special case of witness
problems, and some of them are important for applications to
cryptography: given a property $\mathcal{P}$  and the information
that there are objects with the property $\mathcal{P}$, find
something ``material" establishing the property $\mathcal{P}$; for
example, given two conjugate elements of a group, find a conjugator.
\end{itemize}

All decision problems in group theory have a ``companion"  witness
version, and most of them also have a search  version, and it is the
purpose of this section to illustrate this point by using some of
the most popular algorithmic problems. Below $F$ denotes a free
group with the set (``alphabet") $X$ of free generators, and
$gp_F(R)$ denotes the normal closure of a set $R$ of elements of $F$
in $F$.

\begin{enumerate}

\item Let $G=F/gp_F(R)=<X; R>$ be a finite (or more generally, recursive)
presentation of a group $G$. The already mentioned word (decision)
problem for $G$ is: given a word $w$ in the alphabet $X$, find out
whether or not $w$ is equal to 1 in $G$ or, equivalently, whether or
not $w$ is in the normal closure of $R$.

The word witness problem then is: given that a word $w$ is in the
normal closure of $R$, find a proof (a ``witness") of  that fact.

A particular way of proving it would be to find an expression of $w$
as a product of words of the form $f_i^{-1}r_i^{\pm 1}f_i$, $r_i \in
R$; this can therefore be considered a relevant search problem.

We  note that the  word search problem always has a recursive
solution because one can recursively enumerate all products of
defining relators, their inverses and conjugates. However, the
number of factors in such a product required to represent a word of
length $n$ which is equal to 1 in $G$, can be very large compared to
$n$. If one now considers all   words $w$ of length at most $n$ in
$gp_F(R)$, then the minimum number of conjugates of $r_i^{\pm 1}$
required to express those $w$ gives rise to a function $f(n)$,
termed the {\em isoperimetric function} of the group $G=F/gp_F(R)$.
It provides  one of the possible  measures of complexity of the word
search problem for $G$. It is possible to show that the
isoperimetric function  can be made as complicated a function as one
wishes (see \cite{birget, sapir}). Furthermore, if in a group $G$
the word problem is recursively unsolvable, then the length of a
proof verifying that $w=1$ in $G$ is  not bounded by any recursive
function of the length of $w$.

\item The conjugacy (decision) problem for $G$ is: given two words $w_1,
w_2$, find out whether or not there is a word $g$ such that the
words $g^{-1}w_1g$  and $w_2$ represent the same element of the
group $G$. If they do, then we say that the elements of  $G$
represented by $w_1$  and $w_2$ are {\it conjugate} in $G$.

The conjugacy witness  problem then is: given two words $w_1, w_2$
representing conjugate elements of  $G$, find a  proof (a
``witness") of  the fact that the elements are conjugate.

One of the ways  of proving it would be to find  a particular word
(a conjugator) $g$ such that $g^{-1}w_1g$  and $w_2$ represent the
same element of $G$; this is the conjugacy search problem.

Again, the  conjugacy search problem always has a recursive solution
because one can recursively enumerate all conjugates of a given
element, but as with the word search problem, this kind of solution
can be extremely inefficient.

 We note, in passing, that several cryptographic primitives based on the (alleged)
computational hardness of the conjugacy search problem (in
particular, in braid groups) have been suggested, including
\cite{AAG, GrigShpil, KLCHKP}.

\item The  decomposition (search) problem is:
 given two elements $w_1, w_2$ of a  group  $G$ and two
 subgroups $A, B \le G$ (not necessarily distinct), find  elements
 $x \in A, ~y \in B$ such that $w_1 = x w y$ in $G$,  provided at least one such
pair of elements exists.

 We note that  {\it some} $x$ and $y$   satisfying the equality
$x\cdot g \cdot y =   h$ always exist (e. g.  $x=1, ~y=g^{-1}h$), so
the point is to have them satisfy the conditions $x \in A, ~y \in
B$. We therefore will not usually refer to this problem as a {\it
subgroup-restricted} decomposition search problem because it is
always going to be subgroup-restricted; otherwise it does not make
much sense.

A special case of the decomposition search problem, where $A=B$, is
also known as the {\it double coset problem}.

  The corresponding decision problem is not among  problems
traditionally studied in group theory. The search version (which
generalizes the conjugacy search problem), on the other hand, has
been recently used in several cryptographic protocols including
\cite{KLCHKP, SU3}.

\item Another special case of the decomposition problem is
the {\it factorization  problem}: given  an element $w$ of a group
$G$ and two subgroups $A, B \leq G$, find out whether or not there
are two elements $x \in A$ and $y \in B$ such that $x\cdot y = w$.

The {\it factorization search problem} then is: given  an element
$w$ of a recursively presented group $G$ and two recursively
generated subgroups $A, B \leq G$, find any two elements $x \in A$
and $y \in B$ that would satisfy $x\cdot  y = w$, provided at least
one such pair of elements exists.

\item The subgroup membership (decision) problem is: given a group $G$, a subgroup
$H$ generated by  $h_1, \dots, h_k$, and an element $g \in G$, find
out whether or not $g \in H$.

We note that the membership problem also has a less descriptive
 name, ``the generalized word problem".

The subgroup membership witness   problem then  is: given a group
$G$, a subgroup $H$ generated by  $h_1, \dots, h_k$, and an element
$h \in H$, find a  proof   of  the  fact   that $h \in H$.

An  obvious particular way of proving it would be to find  an
expression of $h$ as  a word in  $h_1, \ldots, h_k$; this is the
subgroup membership    search problem.

\item The  isomorphism (decision)  problem is:
given two finitely presented groups $G_1$ and $G_2$, find out
whether or not they are isomorphic.

The  isomorphism witness   problem  is: given two isomorphic
finitely presented groups $G_1$ and $G_2$, find a proof  of  the
fact   that they are isomorphic.

An  obvious  particular way of proving it would be finding  an
isomorphism between the two groups;  this is the isomorphism search
problem.

\item The  automorphism (decision)  problem is:  given a group $G$ and
two elements $u, v$ of   $G$, find out whether or not there is an
automorphism $\alpha$  of   $G$ such that $\alpha(u)=v$. This is
sometimes also called the automorphic conjugacy problem.

The  automorphism witness   problem  is: given $u, v \in G$, find a
proof  of  the existence of $\alpha \in Aut(G)$ such that
$\alpha(u)=v$, provided at least one such $\alpha$ exists.

Again, an  obvious  particular way of proving it would be finding a
particular automorphism $\alpha$ such that $\alpha(u)=v$; this is
the automorphism search problem.

Long time ago, Whitehead has solved the automorphism decision
problem in a free group $F_r$ of any finite rank $r \ge 2$ (see e.g.
\cite{LS}), and that was one of the most important contributions to
combinatorial group theory in the first half of the 20th century.
But only recently the question about computational complexity of
this problem has been raised \cite{MS_orbits}   and studied
\cite{KSS}, \cite{Khan}, \cite{Lee_2006_1}, \cite{Lee_2006_2}. It is
still unknown, at the time of this writing (cf. \cite[Problems
(F25), (C2)]{BMS}), whether this decision problem is in the class
${\mathbf P}$ (with respect to the lexicographic length of the
inputs) or even ${\mathbf {NP}}$  if $r \ge 3$; it is in the class
${\mathbf P}$ if $r = 2$, according to \cite{Khan} and
\cite{Lee_2006_1}. On the other hand, generically, i.e., on ``most"
inputs, the ``no" answer can be given in linear time, see
\cite{KSS}.

\item The  endomorphism (decision)  problem is: given a group $G$ and
two elements $u, v$ of   $G$, find out whether or not there is an
endomorphism $\alpha$  of   $G$ such that $\alpha(u)=v$.

Relevant witness and search problems are similar to those for the
automorphism   problem.

We point out that the endomorphism problem translates into an {\it
equation} of a special form in the given group $G$. Equations in
groups are a major subject of research, but it is outside of the
scope of the present survey.

\end{enumerate}

Now we make one general observation. Decision problems usually
naturally split into the ``yes" and ``no" parts, and the ``yes" part
of most popular decision problems in group theory  usually  has a
recursive solution; for example, the ``yes" part of the word problem
has a recursive solution because, given a recursive presentation of
a group $G$, the set of all words $w$ such that $w=1$ in $G$ is
recursively enumerable. The same can be said about the ``yes" part
of the conjugacy problem, the isomorphism problem, etc. At the same
time, the ``no" part of these problems is typically {\it not}
recursively enumerable in general. However, one can still ask for a
proof (a ``witness") of the fact that, say, a given word $w$ is not
equal to 1 in $G$, or a given pair of words represent non-conjugate
elements of $G$, etc.  We call the corresponding search problems the
non-identity witness problem (because calling it the ``non-word
witness problem" would be kind of ridiculous) and the non-conjugacy
witness problem, respectively. Similarly, one can consider
non-membership witness problem, non-isomorphism witness problem,
etc.

As we have pointed out before, in general there is no recursive
procedure for enumerating all words $w$ such that $w \ne 1$ in $G$,
or all words representing elements that do not belong to a given
subgroup of $G$, etc. Of course, if, for example, $G$ has solvable
word problem, then enumerating all words $w \ne 1$ in $G$ is
possible by an  obvious procedure. However, what we are looking for
here is a more general way of proving $w \ne 1$ that would be
applicable also to  ``many"  groups with unsolvable word problem.
One fairly general approach to proving $w \ne 1$ in $G$ would be to
exhibit a ``nice" factor group of $G$ (often just the abelianization
$G/[G, G]$ would work) where $w \ne 1$. This approach is discussed
in detail in \cite{KMSS1}, and it also works for the non-conjugacy
witness problem and for the non-membership witness problem. Still,
it would be quite interesting to find other sufficiently general
methods for proving non-identity, non-conjugacy, etc.

At the same time, it would be quite interesting (and useful) to have
a general way (applicable to {\it any} non-trivial  group $G$) of
proving $w \ne 1$ at least for some particular words $w$ (depending
on $G$). This may be regarded as a special case of the
non-isomorphism witness problem, namely,  proving that a given group
is non-trivial:

\begin{problem}\label{non-trivial element} (M. Chiodo \cite{Chiodo:2010}, \cite{BMS})
Is there  a general procedure to produce a non-trivial element from
a finite presentation of a non-trivial group?
\end{problem}

This problem is discussed in \cite{Chiodo:2010}, where a special
case is settled; namely, it is shown that there is no general
procedure to pick a non-trivial generator from a finite presentation
of a non-trivial group. We note here the importance for
cryptographic applications of any progress on Problem
\ref{non-trivial element} in the positive direction.

Building on the same idea, one can also ask:

\begin{problem}\label{not in the subgroup}
Is there  a general procedure to produce an element  that does not
belong to a given (finitely generated) proper subgroup of a given
finitely presented group, provided such elements exist?
\end{problem}

In a somewhat different direction:

\begin{problem}\label{free} Given a
finitely presented group  $G$, elements  $h_1, \ldots, h_k \in G$,
and the information that $h_1, \ldots, h_k$ freely generate a free
subgroup of $G$, find a proof (a ``witness") of that fact.
\end{problem}

It is a matter of taste whether to consider the connotation of the
property alluded to in this problem ``positive" or ``negative". By
``negative" here we mean the absence of nontrivial relations between
$h_1, \ldots, h_k$. We note that if there are relations between
$h_1, \ldots, h_k$, then we can eventually find one by going over
words in $h_1, \ldots, h_k$ and initiating an algorithm for the
``yes" part of the word problem for each.

 To appreciate the difficulty of Problem \ref{free},  the reader may look
at \cite{Collins} to see that even in such well-studied groups as
braid groups, it is not easy to prove that squares of two
``neighbor" braid generators freely  generate a free group. We also
note that, according to \cite{GMO}, a ``random" finite set of
elements of a nonelementary  hyperbolic group $G$ is ``very likely"
to be a set of free generators for a  free subgroup of $G$. This is
consistent with observations made  in \cite{KMSS} concerning the
genericity of the ``no" answer to several other algorithmic problems
in groups, including the word problem, conjugacy problem, etc.

 To conclude this section,  we point out that some specific
problems may provide examples of natural group-theoretic decision
problems with both the ``yes" and ``no" parts  nonrecursive, which
would be of great interest. Here we can offer some candidate
problems of that kind.

\begin{problem}
Is the set of all finitely presented metabelian groups recursively
enumerable?
\end{problem}

A group is called metabelian if its commutator subgroup is abelian.
 The set of finitely presented non-metabelian groups is known to be
nonrecursive, see e.g. \cite{adian}. At the same time, there is no
obvious way to recursively enumerate all finitely presented
metabelian groups because many metabelian groups are not finitely
presented, so it is not clear how to specifically enumerate just the
finitely presented ones.

The relevance of this problem to the present survey is due to the
fact that it may provide a natural  example of a witness problem in
group theory that is algorithmically unsolvable. There are many
algorithmically unsolvable decision problems in group theory (see
e.g. \cite{adian} or \cite{KS-survey}), but the following might be
the first natural example of an unsolvable witness problem:

\begin{problem}
Given a finitely presented group  and the information that it is
metabelian, find a proof (a ``witness") of that fact.
\end{problem}

Another interesting decision  problem that might have both the
``yes" and ``no" parts nonrecursive is: given two finitely presented
groups $G_1$ and  $G_2$, is there an injective homomorphism (an
embedding) of $G_1$ into $G_2$? This problem is known to have a
negative answer, but the point is, again, that it might have both
the ``yes" and ``no" parts nonrecursive, as was suggested to the
author by D. Groves. We note that without the word ``injective", the
``yes" part of this problem would have an affirmative answer, i.e.,
all homomorphisms of $G_1$ into $G_2$ are recursively enumerable.

Thus, we have the following natural  witness problem that may be
algorithmically unsolvable:

\begin{problem}\label{Groves} (D. Groves)
Given two finitely presented groups  $G_1$  and  $G_2$ and the
information that there is an injective homomorphism (an embedding)
of $G_1$ into $G_2$, find a proof (a ``witness") of that fact.
\end{problem}

We note that Chiodo \cite{Chiodo:2010} has recently proved that
there is no algorithm that, on input of finite presentations of two
groups and information that one of them embeds into the other,
outputs an explicit embedding. Therefore, the embedding search
problem is algorithmically unsolvable! This, however, does not
necessarily provide a negative answer to Problem \ref{Groves} above
because there might be other ways to prove the existence of an
embedding (for example, if $G_1$ is a cyclic group of order $n$,
then finding an element of order $n$ in $G_2$    would be such a
proof), but this is a serious argument in favor of a negative answer
nonetheless.


 Another example of a similar kind was reported in the same paper
\cite{Chiodo:2010}: given a finite presentation of a group $G$  and
information that $G$ has an element of a finite order $n \ge 2$,
there is, in general, no algorithm to find a particular element of
order $n$. In fact, it was shown in \cite{Chiodo:2010} that there is
no algorithm to even find {\it any} torsion element in $G$. Again,
this does not necessarily imply that there is no {\it proof} (or
``witness") of the existence of an element of order $n$.

\section{Stratification}
\label{strat}

In this section, we discuss the concept of {\it stratification},
which is important (and also independently interesting) from
theoretical point of view, but at the same time it provides a bridge
between ``more theoretical" class of decision problems and a ``more
practical" class of search problems.

In the end of the previous section, we gave   examples of
algorithmically unsolvable search problems. However, as we have
pointed out also in the previous section, ``standard" search
problems in group theory are algorithmically solvable, so the
question of interest is about the {\it computational complexity} of
search problems.

To allocate a search problem to one of the established complexity
classes (such as  ${\mathbf P}$,  ${\mathbf {NP}}$, etc.), one needs
to convert it to a decision problem. A  standard way of doing it is
to provide some kind of stratification of a possible search outcome
for the search problem at hand. Here is an example of how one can
convert the conjugacy search problem to a decision problem:

\begin{problem}\label{conj_strat}
Given two words $w_1, w_2$ representing conjugate elements of  $G$,
and a positive integer $k$, is there a word $g$ of length at most
$k$ such that $g^{-1}w_1g$ and $w_2$ represent the same element of
$G$?
\end{problem}

Of course,  computational complexity of this problem may depend on
$k$, among other  things. More importantly:

\medskip

\noindent {\bf Warning.} The conjugacy search problem is
algorithmically solvable in any recursively presented group $G$,
whereas Problem \ref{conj_strat} may not be if the word problem in
$G$ is algorithmically unsolvable.

To see that the conjugacy search problem is always solvable, we use
a straightforward algorithm: recursively enumerate all words in the
given generators of $G$, then go over all these words $g$ one at a
time, comparing  $g^{-1}w_1g$ to $w_2$ by using the fact that the
``yes" part of the  word problem is solvable in any recursively
presented group $G$. The crucial point here is   that when we say
``comparing" two elements, we mean {\it initiating} the obvious
procedure for the ``yes" part of the word problem. However, after
initiating such a procedure we do not just sit there waiting for a
result because we do not know how long we have to wait (perhaps
indefinitely); instead, we move on to the next word, initiate the
relevant procedure for the ``yes" part of the word problem, etc.

If, however, we try to use the same procedure for Problem
\ref{conj_strat}, this may not work if the word problem in $G$ is
algorithmically unsolvable. Indeed, suppose we go over all words $g$
of length at most $k$ (in the given generators of $G$) one at a
time, comparing  $g^{-1}w_1g$ to $w_2$ by virtue of  the fact that
the ``yes" part of the  word problem is solvable in $G$. That means,
we have initiated a number (which is, incidentally, exponential in
$k$) of relevant procedures. After that, all we can do is sit there
hoping that one of the initiated procedures would terminate.
However, if the word problem in $G$ is unsolvable, there is no
recursive bound on the run time of any of our procedures, which
means that Problem \ref{conj_strat} is, in general, unsolvable. Note
also that if $k=0$, then Problem \ref{conj_strat}  becomes
equivalent to the word problem in $G$.
\medskip

Thus, the bottom line is: {\it Problem \ref{conj_strat} may not be
algorithmically solvable, while  the conjugacy search problem always
is.} This leaves the problem of allocating the conjugacy search
problem in a given group to one of the complexity classes ``somewhat
open".

\medskip

We note that stratification of a search problem is  often not
unique. Examples of different stratifications of the word search
problem are given in our  Section  \ref{exp_wp}. Examples of
different stratifications of the isomorphism search problem are
given   below.

Here we give another example of a stratification that is relevant to
a ramification of the word problem sometimes called the {\it
geodesic problem}:

\begin{problem}
Given  a word $w$, a  group $G$, and a positive integer $k$, is
there a word $g$ of length at most $k$, which is equal to $w$ in
$G$?
\end{problem}

This problem was shown to be ${\mathbf {NP}}$-hard in some groups
$G$, including, somewhat surprisingly, the free metabelian group of
rank 2 \cite{MRUV}.

Yet another example of a stratification, of a different nature, is
relevant to the isomorphism search problem. There is an obvious
stratification by the sum of the lengths of images of the generators
under a given isomorphism. A much more interesting stratification
however is provided by {\it Tietze transformations}. It is known
that  two groups given by their finite presentations are isomorphic
if and only if one can get from one of the presentations to the
other by a sequence of Tietze transformations; see our subsection
\ref{Tietze} for more details. Therefore, one can stratify the
isomorphism search problem by the length of a sequence of Tietze
transformations establishing an isomorphism between groups.

We emphasize once again that a stratification of a given search
problem is typically not unique, and selecting a ``good" (for
specific purposes) stratification of one search problem or another
can be an important problem of independent interest.

We also note that there is the following important connection
between complexity of a decision problem and that of the associated
witness problem. Suppose the decision problem is: does a given input
$S$ have a property ${\mathcal P}$? If there were an  algorithm
${\mathcal A}$ that would produce, for any $S$ having  property
${\mathcal P}$, a proof of that fact in time bounded by a known
function $f(|S|)$ in the ``size" $|S|$ of $S$, then, given an
arbitrary $S'$, we could run the algorithm ${\mathcal A}$ on $S'$,
and if it would not produce a proof of $S'$ having the property
${\mathcal P}$ after running over the time $f(|S'|)$, we could
conclude that $S'$ does not have the property ${\mathcal P}$,
thereby solving the corresponding decision problem in time
$f(|S'|)$. In particular, a polynomial-time solution of a
witness/search problem implies a polynomial-time solution of the
relevant decision problem.

\subsection{Tietze transformations: elementary isomorphisms}
\label{Tietze}

In this section, we briefly explain a rather nontrivial
stratification of the isomorphism search problem by means  of Tietze
transformations, to illustrate a point that we made above, namely
that a stratification of a search problem may sometimes be rather
nontrivial and may by itself open interesting research avenues.
These are ``elementary isomorphisms": any isomorphism between
finitely presented groups is a composition of Tietze
transformations.

 Tietze introduced isomorphism-preserving elementary transformations that can be applied to
groups presented by generators and relators. They are of the
following  types (here we do not worry about some of the
transformations possibly being redundant).
\begin{description}
\item[(T1)]
\emph{Introducing a~new generator}: Replace $\langle x_1,x_2,\ldots
\mid r_1, r_2,\dots\rangle$ by $\langle y, x_1,x_2,\ldots \mid
ys^{-1}, r_1, r_2,\dots\rangle$, where $s=s(x_1,x_2,\dots )$ is an
arbitrary element in the generators $x_1,x_2,\dots$.
\item[(T2)]
\emph{Canceling a~generator} (this is the converse of (T1)): If one
has a~presentation of the form $\langle y, x_1,x_2,\ldots \mid q,
r_1, r_2,\dots\rangle$, where $q$ is of the form $ys^{-1}$, and $s,
r_1, r_2,\dots $ are in the group generated by $x_1,x_2,\dots$,
replace this presentation by $\langle x_1,x_2,\ldots \mid r_1,
r_2,\dots\rangle$.
\item[(T3)]
\sloppy \emph{Applying an automorphism}: Apply an automorphism of
the free group generated by $x_1,x_2,\dots $ to all the relators
$r_1, r_2,\dots$.
\item[(T4)]
\emph{Changing defining relators}: Replace the set $r_1, r_2,\dots$
of defining relators by another set $r_1', r_2',\dots$ with the same
normal closure. That means, each of  $r_1', r_2',\dots$ should
belong to the normal subgroup generated by $r_1, r_2,\dots$, and
vice versa.
\end{description}

Tietze has proved (see e.g. \cite{LS}) that two groups $\langle
x_1,x_2,\ldots \mid r_1, r_2,\dots\rangle$ and $\langle
x_1,x_2,\ldots \mid s_1, s_2,\dots\rangle$ are isomorphic if and
only if one can get from one of the presentations to the other by a
sequence of transformations \textup{(T1)--(T4)}.

 For each Tietze transformation of the types \textup{(T1)--(T3)}, it is easy to obtain an explicit
 isomorphism (as a mapping on
 generators) and its inverse. For a Tietze transformation of the
 type \textup{(T4)}, the isomorphism is just the identity map.
We would like here to make Tietze transformations of the
 type \textup{(T4)} recursive, because {\it a priori} it is not clear
 how to actually implement them:

\medskip

\noindent {\bf (T4$_1$)} In the set $r_1, r_2,\dots$, replace some
$r_i$ by one of the:  $r_i^{-1}$,  $r_i r_j$, $r_i r_j^{-1}$, $r_j
r_i$, $r_j r_i^{-1}$, $x_k^{-1} r_i x_k$, $x_k r_i x_k^{-1}$, where
$j \ne i$, and $k$ is arbitrary.
\medskip

\noindent {\bf (T4$_2$)} To the set $r_1, r_2,\dots$, add one of the
elements specified in (T4$_1$), without the restriction $j \ne i$.

\medskip

\noindent {\bf (T4$_3$)}  (this is the converse of (T4$_2$)) From
the set $r_1, r_2,\dots$, remove an element if it can be obtained
from other elements as specified in (T4$_2$).

\medskip

It is known \cite{M_AC} that (T4$_1$), (T4$_2$), (T4$_3$) are indeed
sufficient to implement any (T4). We note that finding a sequence of
Tietze transformations between two given presentations of isomorphic
groups is similar to finding a sequence of Andrews-Curtis ``moves"
reducing a balanced presentation of the trivial group to the
``standard" one, which is relevant to the famous Andrews-Curtis
conjecture, well known in combinatorial  group theory and topology.

\section{Computational approach to search problems}
\label{exp_wp}

In this section we explain some ideas, due to Ushakov \cite{U},
behind  ``practical", or computational, approaches to search
problems in group theory. We point out, up front, that we consider
it satisfactory when  a proposed algorithm is efficient on ``most"
inputs, while on a ``negligible" set of inputs it may be inefficient
or may even not terminate. Here ``most" and ``negligible" have
precise meanings, as defined in \cite{KMSS1}.

As explained in the previous sections, even though  decision and
search problems have a lot in common, their computational paradigms
are different. For instance, the most popular search problems (like
the word search problem, the conjugacy search problem) are always
solvable, whereas the corresponding   decision problems may not be.
Theoretical solvability of search problems, however, does not
usually help much in practical implementations   because the upper
bound for the runtime of a search problem algorithm remains a
non-recursive function if the relevant decision problem is
undecidable. On the other hand, since we look at the problems from
the practical point of view, we assume that the instances of the
problem are somehow sampled by some procedure, and the procedure
``knows'' that the sampled instance is a positive instance of the
problem, i.e., it has a proof that the instance is positive. This
spreads out the complexity  ``more evenly'' between two entities,
the one which generates a positive instance of the problem and the
one which finds a proof for that instance.

In summary, we treat a search problem here as a two-party game. We
assume that one party (Alice) generates a positive instance of the
problem, and the other party (Bob) attempts to find a ``witness"
(i.e., a proof) for that instance. In this setting, a natural
analysis of the problem would be the comparison of run times of an
algorithm required for Alice to generate the instance  versus that
for Bob to find  a witness.

The way Alice generates her instances is crucial here. Some
instances can be structurally more complex than other, and hence it
might be more difficult to generate them, i.e., their generation
takes more time. We would like to point out that there is no way to
perform ``uniform" generation process for positive instances in case
the decision problem is undecidable because having  a
``non-recursive time" for generating instances does not make sense.
Therefore, it would be a natural approach to consider the ``size" of
a positive instance of a problem to be the time required to generate
that particular instance. For example:
\begin{itemize}
    \item
(Word Search Problem -- 1) If Alice generates a word $w$ that
represents the identity of $\gp{X; R}$ as a product $w =
\prod_{i=1}^k c_i^{-1} r_i c_i$, then it is natural to say that
$\sum_{i=1}^k (2|c_i| + |r_i|)$ is the size of $w$. This way of
generating   trivial elements of $G$ is natural since it comes from
the definition of the word problem.
    \item
(Word Search Problem -- 2) A slightly different approach to
generating positive instances of the word problem for a given
presentation $\gp{X; R}$ is the following iterative procedure. We
construct a sequence of group words $w_0,\ldots,w_n,w$, where
$\varepsilon=w_0$, $w_i$ is obtained from $w_{i-1}$ by inserting an
element of the form $h_i^{-1} r_i h_i$, and $w$ is obtained by
freely reducing $w_n$. The size of such an instance would be
$\sum_{i=1}^k (2|h_i| + |r_i|)$.
    \item
(Conjugacy Search Problem) If Alice generates a word $v$
representing the conjugate element of a word $u$ of $\gp{X; R}$ as a
product $c^{-1}u c \prod_{i=1}^k c_i^{-1} r_i c_i$, then it is
natural to say that  $2|u|+2|c| + \sum_{i=1}^k (2|c_i| + |r_i|)$ is
the size of the pair $(u,v)$.

  \item
(Membership Search Problem) Let $H$ be a subgroup of $\gp{X; R}$
generated by $\{h_1,\ldots,h_k\} \subset F(X)$. Alice can generate
elements of $H$ as follows. She constructs a sequence of words
$w_0,\ldots,w_n, w$, where $w_0 = \varepsilon$; every $w_{i+1}$ is
obtained either by multiplying $w_i$ by some $h_{j_i}$ on the left
or on the right, or by inserting  a word of the form $c_i^{-1} r_i
c_i$ inside $w_i$; and $w$ is obtained by reducing $w_n$. It is
natural to say that the size of the   word $w$, which represents an
element of $H$, is
    $$\sum_{i=1}^k
    \left\{
    \begin{array}{ll}
        |h_{j_i}|, & \mbox{if $h_{i_j}$ is inserted at the $i$th step};\\
        2|c_i|+|r_i|, & \mbox{if $c_i^{-1} r_i c_i$ is inserted at the $i$th step}.\\
    \end{array}
    \right.
    $$

\end{itemize}

Slightly modified methods were analyzed in  \cite{U}.



\end{document}